\documentclass[12pts]{article}
\usepackage{amsfonts}
\usepackage{amsmath}
\usepackage{amssymb}
\usepackage{graphicx}%
\providecommand{\U}[1]{\protect\rule{.1in}{.1in}}

\begin{document}

\title{Remarks on the Choquet Integral Calculus on $[a, t]$, with $a\in \mathbb{R}$}
\author{Sorin G. Gal \\
Department of Mathematics and Computer Science, \\
University of Oradea, \\
Universitatii 1, 410087, Oradea, Romania\\
E-mail: \textit{galsorin23@gmail.com, galso@uoradea.ro}}

\maketitle

\begin{abstract}
In this note, we extend the considerations for the Choquet integral calculus on the interval $[0, t]$ introduced in \cite{Su}, \cite{Su3}, to the case of an interval $[a, t]$, with arbitrary $a\in \mathbb{R}$.
\end{abstract}



\section{Introduction}

In two very seminal papers \cite{Su}, \cite{Su3}, the basis of a theory concerning the calculations of the continuous Choquet integral
and the inverse problem of the Choquet integral were posed, in both cases considered with respect to distorted Lebesgue measures on the nonnegative real line. Then, in \cite{Rida} the theory was extended to the case of non-monotonous nonnegative functions.

Let us shortly recall the main elements of this theory.

For ${\cal{B}}$  the
smallest $\sigma$-algebra including all the closed intervals in $\mathbb{R}_{+}=[0, +\infty)$, let $\mu:{\cal{B}}\to \mathbb{R}_{+}$ be a (non-null) capacity. Denote simply by ${\cal{F}}^{+}$ the space of all nonnegative, nondecreasing and continuous functions defined on  $[0, +\infty)$.

Starting from the Choquet integral equation
\begin{equation}\label{eq22}
f(t)=f(0)+(C)\int_{0}^{t}g(s)d \mu(s), t\in [0, +\infty),\, t\ge 0,
\end{equation}
in \cite{Su}, \cite{Su3} were posed the following  three problems :

1) Calculation of the Choquet integral: calculate $f\in {\cal{F}}^{+}$
for given $g\in {\cal{F}}^{+}$ and $\mu$ ; without any loss of generality, in this case we suppose in (\ref{eq22}) that $f(0)=0$.

2) Solution of the Choquet integral equation :
find $g\in {\cal{F}}^{+}$ for given $f\in {\cal{F}}^{+}$ and $\mu$. If exists uniquely, the function $g$ is called the derivative of $f$ with respect to $\mu$ and it will be denoted
with $g=\frac{d f}{d \mu}=f^{\prime}_{\mu}$;

3) identification of the fuzzy measure: identify
$\mu$ for given $f$ and $g$ in ${\cal{F}}^{+}$.

Solutions to the above problems were obtained in \cite{Su}, \cite{Su3}, \cite{Rida}, mainly for distorted Lebesgue measures, with $m$ continuously differentiable on $\mathbb{R}_{+}$ and for continuously  differentiable $f$ and $g$.

The main goal of the present note is to deliver answers to the above three problems.

\section{Statements and Solutions of the Problems}

For $a\in \mathbb{R}$, denote by ${\cal{F}}^{+}[a, +\infty)$ the space of all nonnegative, nondecreasing and continuous functions defined on  $[a, +\infty)$. Note that for simplicity, ${\cal{F}}^{+}[0, +\infty)$ was denoted by ${\cal{F}}^{+}$.

By analogy with the above considerations, starting from the Choquet integral equation
\begin{equation}\label{eq221}
f(t)=f(a)+(C)\int_{a}^{t}g(s)d \mu(s), t\in [a, +\infty),
\end{equation}
where $f\in {\cal{F}}^{+}[a, +\infty)$, we can state the same three problems as in the case when $a=0$.

{\bf Problem 1.} Calculation of the Choquet integral: calculate $f\in {\cal{F}}^{+}[a, +\infty)$
for given $g\in {\cal{F}}^{+}[a, +\infty)$ and $\mu$ ; without any loss of generality, in this case we suppose in (\ref{eq221}) that $f(a)=0$.

An answer to the Problem 1, is the following.

{\bf Theorem 1.1} {\it Suppose that $\mu([a, t])$ is differentiable as function of $t\ge a$ and  $g\in {\cal{F}}^{+}[a, +\infty)$.

(i) We have that
$$f(t)=(C)\int_{a}^{t}g(\tau)d \mu(\tau)=-\int_{a}^{t}\mu^{\prime}_{\tau}([\tau, t])g(\tau)d \tau, \mbox{ for all } t\ge a.$$

(ii) Let $\mu=m \circ \lambda$, where $\lambda$ is the Lebesgue measure and $m$ is differentiable on $\mathbb{R}_{+}$. Then we have
$$f(t)=(C)\int_{a}^{t}g(\tau)d \mu(\tau)=\int_{a}^{t}m^{\prime}(t-\tau)g(\tau)d \tau, \mbox{ for all } t\ge a.$$}

{\bf Proof.} The proofs of (i) and (ii) go exactly as in Proposition 1 in \cite{Su3}. $\hfill \square$

{\bf Remark 1.2.} Concerning the Choquet integral equation on $[0, t]$ in (\ref{eq22}), it is important is it has the hereditary property. With other words, supposing that  $f, g\in {\cal{F}}^{+}$ satisfy the Choquet integral equations on [0, t]
in (\ref{eq22}) (with $f(0)=0$), it is a natural but important question to ask if they also satisfy the Choquet integral equation on $[a, t]$, with arbitrary $a>0$, that is
\begin{equation}\label{eq33}
f(t)=f(a)+(C)\int_{a}^{t}g(s)d \mu(s), t\in [a, +\infty).
\end{equation}
If $f$, $g$ and $\mu$ satisfy the Choquet integral equations in Theorem 1.1, (i), (ii), then the answer to the question is positive. For example,
in the case of Theorem 1.1, (ii) (the proofs for the case (i) is similar)  we get for all $t\ge a$
$$\int_{a}^{t}m^{\prime}(t-\tau)g(\tau)d \tau =\int_{0}^{t}m^{\prime}(t-\tau)g(\tau)d \tau-\int_{0}^{a}m^{\prime}(t-\tau)g(\tau)d \tau=f(t)-f(a),$$
which implies
$$f(t)=f(a)+\int_{a}^{t}m^{\prime}(t-\tau)g(\tau)d \tau, t\in [a, +\infty).$$
But a direct calculation (similar to those in \cite{Su3}, Proposition 1), easily leads to the equality $(C)\int_{a}^{t}g(\tau)d \mu(\tau)=\int_{a}^{t}m^{\prime}(t-\tau)g(\tau)$, for all $t\ge a$, which proves our assertion.

{\bf Remark 1.3.} Another question which naturally arises  is how we could apply the Laplace's transform method for concrete calculation in Theorem 1.1. For this purpose, we will use the property that we can change the variable under the Riemann integral. Thus, denoting $f_{a}(r)=f(r+a)$, $g_{a}(r)=g(r+a)$, $r\ge 0$, obviously that $f_{a}$, $g_{a}$ keep their monotonicity and in fact $f_{a}, g_{a}\in {\cal{F}}^{+}$.
Then, replacing $t\ge a$ by $r+a$ with $r\ge 0$  in both integral equations in Theorem 1, by the change of variable $\tau=\eta+a$, we obtain
$$f_{a}(r)=-\int_{a}^{r+a}\mu^{\prime}_{\tau}([\tau, r+a])g(\tau)d \tau=-\int_{0}^{r}\mu^{\prime}_{\eta}([\eta+a, r+a])g_{a}(\eta)d \eta, \, r\ge 0,$$
and
$$f_{a}(r)=\int_{a}^{r+a}m^{\prime}(r+a-\tau) g(\tau)d \tau=\int_{0}^{r}m^{\prime}(r+a-(\eta+a))g_{a}(\eta)d \eta=\int_{0}^{r}m^{\prime}(r-\eta)g_{a}(\eta)d \eta, \, r\ge 0.$$
Therefore, we have reduced the integral equations in Theorem 1.1, to the forms in the papers \cite{Su}, \cite{Su3}. Then, applying formally to the last integral equation the Laplace transform ${\cal{L}}$, exactly as in \cite{Su}, \cite{Su3} we obtain
\begin{equation}\label{E1}
F_{a}(s)=s M(s) G_{a}(s), s\ge 0,
\end{equation}
where $F_{a}(s)={\cal{L}}(f_{a}(t))(s)$, $G_{a}(s)={\cal{L}}(g_{a}(t))(s)$, $M(s)={\cal{L}}(m(t))(s)$.

This implies $f(t+a)=f_{a}(t)={\cal{L}}^{-1}[s M(s)G_{a}(s)](t)$, $t\ge 0$, and therefore
\begin{equation}\label{E2}
f(t)={\cal{L}}^{-1}[s M(s)G_{a}(s)](t-a), \mbox{ for all } t\ge a.
\end{equation}

{\bf Example 1.4.} In the case of Problem 1, as an example, for $a\in \mathbb{R}$, choose $g(\tau)=\sqrt{\tau-a}$, $\tau\ge a$ and $m(t)=t^{2}/2$. We get $g_{a}(\tau)=\sqrt{\tau}$ and by using a symbolic Laplace transform calculator (https://www.symbolab.com/solver/laplace-calculator/laplace), it follows $G_{a}(s)={\cal{L}}(\sqrt{t})(s)=\frac{\sqrt{\pi}}{2 s^{3/2}}$, $M(s)={\cal{L}}(t^{2}/2)(s)=\frac{1}{s^{3}}$. By using formula (\ref{E1}), it easily follows $F_{a}(s)=\frac{\sqrt{\pi}}{2}\cdot \frac{1}{s^{7/2}}$ and applying now the symbolic inverse Laplace transform (at the same link as above), we obtain
$$f(t+a)={\cal{L}}^{-1}\left (\frac{\sqrt{\pi}}{2}\cdot \frac{1}{s^{7/2}}\right )(t)=\frac{\sqrt{\pi}}{2}\cdot \frac{8 t^{5/2}}{15\sqrt{\pi}}=\frac{4}{15}\cdot t^{5/2}.$$
By (\ref{E2}) it obviously follows $f(t)=\frac{4}{15}\cdot (t-a)^{5/2}$, which is nondecreasing as function of $t\ge a$.

{\bf Problem 2.}  Solution of the Choquet integral equation (\ref{eq221}) :
find $g\in {\cal{F}}^{+}[a, +\infty)$ for given $\mu$ and $f\in {\cal{F}}^{+}[a, +\infty)$ with $f(a)=0$.

Keeping the notations for $F_{a}$ and $G_{a}$ in the previous Remark 1.3 and using the formula (\ref{E1}), we get $G_{a}(s)=\frac{F_{a}(s)}{s M(s)}$, which implies
$g_{a}(t)={\cal{L}}^{-1}\left [\frac{F_{a}(s)}{s M(s)}\right ](t)$.  Therefore,
\begin{equation}\label{E3}
g(t)={\cal{L}}^{-1}\left [\frac{F_{a}(s)}{s M(s)}\right ](t-a), \mbox{ for all } t\ge a,
\end{equation}
if $g$ is found to be in ${\cal{F}}^{+}[a, +\infty)$. In this case, we call $g(t)$ as the derivative of $f(t)$ for $t\ge a$ with respect to $\mu$ and it is denoted with $g(t)=\frac{d f}{d \mu}(t)=f^{\prime}_{\mu}(t)$, $t\ge a$  ;

{\bf Example 1.5.} Take $f(t)=\sqrt{t-a}$ and $m(t)=t^{2}/2$. It follows that $f_{a}(t)=\sqrt{t}$, $F_{a}(s)=\frac{\sqrt{\pi}}{2}\cdot \frac{1}{s^{3/2}}$,
$M(s)=\frac{1}{s^{3}}$, $s M(s)=\frac{1}{s^{2}}$,
$$g(t+a)=g_{a}(t)={\cal{L}}^{-1}\left (\frac{F_{a}(s)}{s M(s)}\right )(t)=\frac{\sqrt{\pi}}{2}{\cal{L}}^{-1}(\sqrt{s})(t)=\frac{\pi}{4}\cdot \frac{1}{t^{3/2}},$$
which leads to $g(t)=\frac{\pi}{4}\cdot \frac{1}{(t-a)^{3/2}}$. Since $g(t)$ is nonincreasing and not nondecreasing as required, it follows that dose not exist the Choquet derivative on $[a, +\infty)$ of $f(t)=\sqrt{t-a}$ with respect to the set function $\mu(A)=\lambda^{2}(A)/2$, where $\lambda$ is the Lebesgue measure.

Now, if we take, for example,  $f(t)=(t-a)^{3}\sqrt{t-a}$, $t\ge a$ and $m(t)=t^{2}/2$, it follows that $f_{a}(t)=t^{3}\sqrt{t}=t^{7/2}$, $F_{a}(s)=\frac{105 \sqrt{\pi}}{16 s^{9/2}}$, $M(s)=\frac{1}{s^{3}}$, $s M(s)=\frac{1}{s^{2}}$,
$$g(t+a)=g_{a}(t)={\cal{L}}^{-1}\left (\frac{F_{a}(s)}{s M(s)}\right )(t)=\frac{\sqrt{105 \pi}}{16}{\cal{L}}^{-1}\left (\frac{1}{s^{5/2}}\right )(t)=\frac{105 \sqrt{\pi}}{16}\cdot \frac{4}{3 \sqrt{\pi}}\cdot t^{3/2}=\frac{35}{4}\cdot t^{3/2}.$$
This leads to $g(t)=\frac{35}{4}\cdot (t-a)\sqrt{t-a}$, which being nondecreasing as function of $t\ge a$, therefore implies that  $$f^{\prime}_{\mu}(t)=\frac{d f}{d \mu}(t)=\frac{35}{4}\cdot (t-a)\sqrt{t-a}, \mbox{ for all } t\ge a.$$
Recall here that $\mu=m\circ \lambda$, with $\lambda$ the Lebesgue measure.

{\bf Problem 3.} Identify $\mu$ for given $f\in {\cal{F}}^{+}[a, +\infty)$ and $g\in {\cal{F}}^{+}[a, +\infty)$.

In this case, again from (\ref{E1}) we will get
$$m(t)={\cal{L}}^{-1}\left [\frac{F_{a}(s)}{s G_{a}(s)}\right ](t), t\ge a,$$
which gives the solution if $m(t)$ is nonnegative and strictly increasing.

{\bf Example 1.6.} Take $f(t)=(t-a)^{11/2}$, $g(t)=(t-a)^{1/2}$. By using the link mentioned at Example 1.4, we obtain
$F_{a}(s)=\frac{10395 \sqrt{\pi}}{64 s^{15/2}}$, $G_{a}(s)=\frac{\sqrt{\pi}}{2}\cdot \frac{1}{s^{3/2}}$,
$\frac{F_{a}(s)}{s G_{a}(s)}=\frac{10395}{32}\cdot \frac{1}{s^{7}}$ and therefore
$$m(t)=\frac{10395}{32}\cdot {\cal{L}}^{-1}(1/s^{7})=\frac{10395}{32}\cdot \frac{t^{6}}{720}=\frac{1159}{2560}\cdot t^{6},$$
which is strictly increasing.

We end this section with the following important comment.

{\bf Remark 1.7.} If the capacity $\mu$ is invariant at translations (that is $\mu(A_d)=\mu(A)$ for all $d>0$, where $A_{d}=\{d+x ; x\in A\}$), then
the integral equation on $[a, t]$ with $f, g\in {\cal{F}}^{+}[a, +\infty)$
\begin{equation}\label{F_1}
f(t)=(C)\int_{a}^{t}g(\tau)d \mu(\tau), t\ge a,
\end{equation}
can be reduced to the integral equation  $[0, t]$ with $f_{a}(t)=f(t+a)$, $g_{a}(t)=g(t+a)$, $f_{a}, g_{a}\in {\cal{F}}^{+}$,
\begin{equation}\label{F_2}
f_{a}(t)=(C)\int_{0}^{t}g_{a}(\tau) d\mu(\tau), t\ge 0.
\end{equation}
Indeed, replacing $t\ge a$ by $r+a$ with $r\ge 0$ in equation (\ref{F_1}), we obtain
$$f(r+a)=(C)\int_{a}^{r+a}g(\tau)\mu(\tau)=\int_{0}^{\infty}\mu(\{\tau\in [a, t]; g(\tau)\ge \alpha\})d \alpha
=\int_{0}^{\infty}\mu(\{\tau\in [a, r+a]; g(\tau)\ge \alpha\})d \alpha$$
$$=\int_{0}^{\infty}\mu(\{\eta\in [0, r] ; g(\eta+a)\ge \alpha\})d\alpha=(C)\int_{0}^{r}g_{a}(\eta)d \mu(\eta),$$
since $\mu$ is invariant at translation.

Note that all the distorted Lebesgue measures are invariant at translations due to the invariance at translations of the Lebesgue measure.

It is clear that if $\mu$ is not invariant at translations, then in general, the integral equation (\ref{F_1}) cannot be reduced to the integral equation
(\ref{F_2}).


\begin{thebibliography}{99}



\bibitem{Rida}  Ridaoui, M., Grabisch, M.: Choquet integral calculus on a continuous support and its applications, \emph{Oper. Res. Decis.}, {\bf 26} (2016), no. 1, 73-93.

\bibitem{Su} Sugeno, M. : A way to Choquet calculus, \emph{IEEE Trans. Fuzzy Systems}, {\bf 23} (2015) no. 5, 1439-1457.


\bibitem{Su3} Sugeno, M. : A note on derivatives of functions with respect  to fuzzy measures, \emph{Fuzzy Sets and Systems}, {\bf 222} (2013), 1-17.






\end{thebibliography}
\end{document}